\documentclass{article}

\usepackage{arxiv}

\usepackage[utf8]{inputenc} 
\usepackage[T1]{fontenc}    
\usepackage{hyperref}       
\usepackage{url}            
\usepackage{booktabs}       
\usepackage{amsfonts}       
\usepackage{nicefrac}       
\usepackage{microtype}      
\usepackage{lipsum}		
\usepackage{graphicx}
\usepackage{doi}
\usepackage{soul} 

\usepackage[backend=biber,style=numeric]{biblatex}
\usepackage{amssymb}
\usepackage{amsthm}
\usepackage{amsmath}
\usepackage{float}
\usepackage{subcaption}

\DeclareMathOperator*{\argmin}{argmin}
\usepackage{graphicx} 
\usepackage{subcaption} 
\usepackage[linesnumbered]{algorithm2e}
\RestyleAlgo{ruled}
\usepackage{hyperref}
\usepackage{xcolor}

\title{Solving bilevel problems with products of upper- and lower-level variables}

\author{ \href{https://orcid.org/0009-0000-1455-9890}{\includegraphics[scale=0.06]{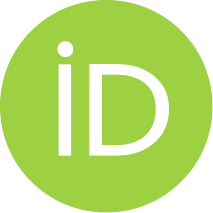}\hspace{1mm}Sina Hajikazemi}
		\\
	Energy Information Networks and Systems\\
	Technical University of Darmstadt\\
	Darmstadt 64283, Germany \\
	\texttt{sina.hajikazemi@eins.tu-darmstadt.de} \\
	\And
	\href{https://orcid.org/0000-0003-3012-9991}{\includegraphics[scale=0.06]{orcid.pdf}\hspace{1mm}Florian Steinke} \\
	Energy Information Networks and Systems\\
	Technical University of Darmstadt\\
	Darmstadt 64283, Germany \\
	\texttt{florian.steinke@eins.tu-darmstadt.de} \\
}

\renewcommand{\headeright}{Preprint}
\renewcommand{\undertitle}{Preprint}

\hypersetup{
pdftitle={Solving bilevel problems with products of upper- and lower-level variables},
pdfsubject={q-bio.NC, q-bio.QM},
pdfauthor={Sina Hajikazemi, Florian Steinke},
pdfkeywords={First keyword, Second keyword, More},
}

\begin{document}
\maketitle

\begin{abstract}
Bilevel programming problems frequently arise in real-world applications across various fields, including transportation, economics, energy markets and healthcare.
These problems have been proven to be NP-hard even in the simplest form with linear upper and lower-level problems.
This paper addresses a specific type of bilevel programming problem where the upper-level is linear, and the lower level includes bilinear terms involving product of variables from both levels.
We propose a new iterative algorithm that addresses this specific class of bilevel problems by penalizing the duality gap and linearizing the bilinear terms.
The effectiveness of the algorithm is argued and demonstrated through a numerical example.

\keywords{bilevel programming   \and bilinear \and penalty \and linearization}
\end{abstract}

\section{Introduction}
\label{sec:intro}
Bilevel programming problems frequently arise in real-world applications across various fields, including transportation, economics, energy markets and healthcare.
This problem is NP-hard even in the simplest form with linear upper and lower-level problems\cite{Hansen1992}.

A common approach to solve the LP-LP bilevel programming problem to the global optimum is to replace the lower-level optimization problem with its KKT conditions, which turns it into a single-level optimization problem.
However, the complementary slackness constraints involve nonlinear terms.
These constraints can be linearized using either binary variables and a big-M value \cite{fortuny1981representation} or SOS1 constraints \cite{kleinert2023sos1}.
Both approaches are typically very slow to solve and are not suitable for large problems.
In addition, finding an appropriate big-M value is very challenging;
\cite{kleinert2020there} shows that finding an appropriate big-M value is as hard as solving the original problem.
Another method to solve this problem is the Penalty Alternating Direction Method (PADM) \cite{kleinert2021computing}. Unlike the previous methods
it only converges to a partial optimum, but it is computationally feasible for large problems and is shown to yield solutions close to the global optimum in a large fraction of sample problems. 

The problem studied in this research, is linear in the upper level and nonlinear in the lower level, with bilinear terms of the upper and lower-level variables in the constraints.
It is easy to prove that any LP-LP bilevel problem can be reduced to this problem by defining and fixing some lower level variables. Thus, the studied problem is harder than an LP-LP bilevel problem.
To the best of our knowledge, there is no prior research that addresses this specific class of problems. 

To solve the proposed, computationally hard problem,
we design a scalable, iterative algorithm called Penalty Adaptive Linearization Method (PALM).
This iterative algorithm penalizes the optimality gap of the lower-level problem 
similar to PADM \cite{kleinert2021computing}, balancing the optimality of the upper- and lower-level problems with increasing weight on the lower-level optimality gap during the course of the algorithm execution.
During this process, it also linearizes the bilinear terms in the constraints and iteratively readjusts the linearization.

The rest of the paper is organized as follows.
Section \ref{section:formulation} formulates the problem.
Section \ref{section:method} introduces the proposed algorithm.
Section \ref{section:example} demonstrates the algorithm with a small numerical example.
Finally, section \ref{section:conclusion} discusses the results and outlines possible future work.

\section{Problem Formulation}\label{section:formulation}
The studied problem consists of two levels.
The upper-level problem is given in (\ref{upper}), where $A \in \mathbb{R}^{p\times mn}$, $B \in \mathbb{R}^{p\times n}$, $c \in \mathbb{R}^{mn}$, $d \in \mathbb{R}^{n}$ and $a \in \mathbb{R}^{p}$ are the coefficients.
The variables $X \in \mathbb{R}^{m \times n}$ and $y \in \mathbb{R}^{n}$ correspond to the upper and lower-level problems, respectively.
Additionally, function $vec$ is used to vectorize the elements of a matrix,
and $S(X)$ denotes the point to set mapping from $X$ to the optimal solutions of the lower-level problem.

\begin{subequations}\label{upper}
    \begin{align}
    \min_{X,y} \quad  & c^T vec(X) + d^Ty \label{obj-upper-level} \\
    & A\, vec(X) + B y\; 
     \geq a \label{upper-constr-connecting} \\
    & y \in S(X) \label{upper-constr-lower-opt}
    \end{align}
\end{subequations}

The lower-level problem is given in (\ref{primal}), where $C \in \mathbb{R}^{m\times n}$, $b \in \mathbb{R}^{m}$ and $e \in \mathbb{R}^{n}$ are the coefficients.
Constraint (\ref{upper-constr-lower-opt}) in the upper-level problem ensures that $y$ attains the optimal value of the lower-level problem.
\begin{subequations}\label{primal}
    \begin{align}
    \min_{y}\quad   &e^T y \\
    s.t.\quad & (C+X)\, y \geq b \label{primal-constr-1}
    \end{align}
\end{subequations}

\section{Proposed Solution Method}\label{section:method}
To obtain a single-level formulation of the bilevel problem, the lower-level problem can be replaced by a set of optimality conditions.
For a linear programming problem, these can be primal feasibility, dual feasibility, and strong duality constraints. The strong duality constraint forces the objective function of the primal minimization problem to be less than or equal to the objective function of the dual problem.

Since the lower-level problem (\ref{primal}) is linear with respect to the primal variables $y$, its dual is given by (\ref{dual}), where \(\lambda \in \mathbb{R}^{m}\) denotes the dual variables.
\begin{subequations}\label{dual}
\begin{align}
\max_{\lambda} \quad   &b^T \lambda \\
s.t.\quad & (C+X)^T \lambda = e \label{dual-constr-1}\\
& \lambda \geq 0\label{dual-constr-2}
\end{align}
\end{subequations}
Replacing (\ref{upper-constr-lower-opt}) with the above mentioned optimality conditions yields the following single-level reformulation of (\ref{upper}): 
\begin{subequations}\label{U-model}
\begin{align}
\min_{X,y,\lambda} \quad  & c^T vec(X) + d^Ty \label{U-obj} && \\
s.t. \quad  & (\ref{upper-constr-connecting}) &&\textit{upper-level constraints}  \nonumber\\
& (\ref{primal-constr-1})   &&\textit{primal constraints} \nonumber  \\
& (\ref{dual-constr-1})-(\ref{dual-constr-2}) &&\textit{dual constraints} \nonumber\\
&  e^Ty \leq b^T \lambda &&\textit{strong duality constraint} \label{U-strong-duality-constr}
\end{align}
\end{subequations}
Constraint (\ref{U-strong-duality-constr}) is the strong duality constraint.
Including this constraint along with the primal and dual feasibility constraints strictly enforces the optimality of the lower-level problem.
Problem (\ref{U-model}) is a non-convex quadratically constrained problem which are NP-hard to solve globally in general \cite{sahni1974computationally}. 

Even finding a feasible point is difficult due to the strictness of the strong duality constraint (\ref{U-strong-duality-constr}). 
Hence, we proceed iteratively and  move the strong duality constraint into the objective function as a penalty term, gradually increasing its weight to progressively enforce the optimality of the lower-level problem, similar to the PADM method.
Shifting the strong duality constraint to the objective function as a penalty term with weight $\mu$ results in the following formulation:
\begin{subequations}\label{U-model-penalty}
	\begin{align}
	\min_{X,y,\lambda} \quad  & c^T vec(X) + d^Ty + \mu \left[ e^Ty - b^T \lambda \right] \\
	s.t. \quad  & (\ref{upper-constr-connecting}) &&\textit{upper-level constraints}  \nonumber\\
	& (\ref{primal-constr-1})   &&\textit{primal constraints} \nonumber  \\
	& (\ref{dual-constr-1})-(\ref{dual-constr-2}) &&\textit{dual constraints} \nonumber
	\end{align}
\end{subequations}
Constraints (\ref{primal-constr-1})  and (\ref{dual-constr-1}) both contain bilinear terms which makes problem (\ref{U-model-penalty}) non-convex and difficult to solve.
Our approach is to linearize the bilinear terms in (\ref{U-model-penalty}) and solve the problem iteratively.

To linearize the $Xy$ term in (\ref{primal-constr-1}), we approximate it around a point $(\bar{X},\bar{y})$.
Let $dX \in \mathbb{R}^{m\times n}$ and $dy \in \mathbb{R}^n$ be defined as $dX = X - \bar{X}$ and $dy = y -\bar{y} $.
Then the approximation is as follows:
\begin{align}
	X y  &= \bar{X}\, \bar{y} + \bar{X} \, dy + dX \, \bar{y} + dX \, dy \nonumber \\
	         &\approx \bar{X} y +  dX \, \bar{y} \label{x-y-approx}
\end{align}
Assuming that $dX$ and $dy$ are small, the term $dX \, dy$ is ignored in the second step. 
This approximation helps to linearize the bilinear terms in (\ref{primal-constr-1}).
The same approximation can be applied to linearize the $X^T\lambda$ terms in (\ref{dual-constr-1}).

Replacing the bilinear terms with the linear approximation changes the primal form of the lower-level problem (\ref{primal}) to (\ref{primal-taylor}) and the dual form of the lower-level problem (\ref{dual}) to (\ref{dual-taylor}).

\begin{subequations}\label{primal-taylor}
	\begin{align}
	\min_{y}\quad  & e^T y \\
	& (C+\bar{X}) y + dX \, \bar{y} \geq b   \label{primal-fixed-constr-1}
	\end{align}
\end{subequations}

\begin{subequations}\label{dual-taylor}
	\begin{align}
	 \max_{\lambda} \quad  & b^T \lambda \\
	& (C+\bar{X})^T \lambda   + dX^T \, \bar{\lambda}= e \label{dual-fixed-constr-1}\\
	& \lambda \geq 0 \label{dual-fixed-constr-2}
	\end{align}
\end{subequations}

Substituting the primal and dual feasibility constraints (\ref{primal-constr-1}) and (\ref{dual-constr-1}) in problem (\ref{U-model-penalty}) with (\ref{primal-fixed-constr-1}) and (\ref{dual-fixed-constr-1}),
yields the following formulation:
\begin{subequations}\label{U-model-fixed-penalty}
\begin{align}
\min_{dX,y,\lambda} \quad  & c^T \textit{vec}(\bar{X}+dX) + d^Ty + \mu \left[ e^Ty - b^T \lambda \right]  \label{U-fixed-penalty-obj}\\
s.t. \quad & (\ref{upper-constr-connecting})  \nonumber  \\
& (\ref{primal-fixed-constr-1})   \nonumber \\
& (\ref{dual-fixed-constr-1})-(\ref{dual-fixed-constr-2})  \nonumber
\end{align}
\end{subequations}
The problem given in (\ref{U-model-fixed-penalty}) is linear. Let $S^L(\bar{X},\bar{y},\bar{\lambda})$ denote the set of optimal values of (\ref{U-model-fixed-penalty}).
\bigskip

The proposed algorithm is presented in \ref{alg:ex-padm} and consists of two loops.
The first (outer) loop (starting at line \ref{alg:opt-loop}) evaluates the optimality of the lower-level problem. 
At the end of each iteration, the penalty coefficient $\mu$ is doubled, thereby increasing the priority of the optimality of the lower-level problem over the upper-level objective function.
This loop terminates when the optimality of the lower-level problem is satisfied, that means a feasible solution to the upper-level problem is found.
If it does not terminate after a certain number of iterations, it means that either the problem is infeasible or the algorithm is not able to find a feasible solution.
The second loop (inner) (starting at line \ref{alg:apx-loop}), solves the linear problem \ref{U-model-fixed-penalty}, updates the values of $\bar{X}$, $\bar{y}$ and $\bar{\lambda}$ and reapeats until these values converge.

Fixing the variable $X$ to $\bar{X}$ in problem (\ref{U-model-penalty}) results in a linear problem. Moreover, the constraints involving $y$ and $\lambda$ in are separable in (\ref{U-model-penalty}), so the problem can be decomposed into two subproblems (\ref{sub1}) and (\ref{sub2}), which are solved independently.
\begin{align}\label{sub1}
	\min_{y} \quad  &  d^Ty + \mu  e^Ty   \\
	s.t. \quad & A\, vec(\bar{X}) + B y\; \geq a  \nonumber  \\
	& (C+\bar{X})\, y \geq b   \nonumber
\end{align}
\begin{align}\label{sub2}
	\max_{\lambda} \quad  &  b^T \lambda \\
	s.t. \quad &  (C+\bar{X})^T \lambda = e \nonumber \\
	& \lambda \geq 0 \nonumber
\end{align}

So for a given $\bar{X}$, the solutions to the subproblems (\ref{sub1}) and (\ref{sub2}) determine the values of $\bar{y}$ and $\bar{\lambda}$, respectively.
Let $S^P(\bar{X})$ and $S^D(\bar{X})$ be the optimal solution sets of the subproblems (\ref{sub1}) and (\ref{sub2}), respectively.
If the solutions of (\ref{sub1}) and (\ref{sub2}) are not unique for all $\bar{X}$, then $\bar{y}$ and $\bar{\lambda}$ may jump between different alternative solutions within the sets $S^P(\bar{X})$ and $S^D(\bar{X})$.
This can cause the algorithm to fail to converge in the inner loop. To mitigate this, we select the optimal solution that is closest to the previous one. 
To find it we first compute the optimal objective value by solving the LP. We then resolve it with the objective value fixed to the optimal value and minimizing the distance to the prior solution. Depending on whether the 1-norm or 2-norm is used, this results in a linear or convex quadratic programming problem, respectively. The values of $\bar{y}$ and $\bar{\lambda}$ are then updated according to lines \ref{alg:update-y} and \ref{alg:update-lambda}.
Similarly, when choosing among the alternative optimal solutions in $S^L(\bar{X},\bar{y},\bar{\lambda})$, we choose the solution that results in the smallest change in $\bar{X}$.

The algorithm needs an initial value $X_0$ for the upper-level variable, for which a feasible primal and dual solution exists to start with.
Let $\Omega$ be the set of feasible solutions of problem (\ref{U-model-penalty}).
Then the initial values for $y$ and $\lambda$ are computed as in line \ref{alg:init}.
Note that this problem is linear since $X$ is fixed to $\bar{X}$.

\begin{algorithm}[!h]
\SetKwInOut{Parameter}{parameter}
\caption{Penalty Adaptive Linearization Method (PALM)}\label{alg:ex-padm}
\Parameter{$\mu_0, X_0, \epsilon_{opt}, \epsilon_{apx}$}
$\mu \gets \mu_0$\;
$\bar{X} \gets X_0 $\; 
$\bar{y},\bar{\lambda} \gets \argmin_{y,\lambda}\left\{ e^Ty - b^T \lambda :(\bar{X},y,\lambda) \in \Omega \right\}$ \;\label{alg:init} 
$i\gets 0$\;
\While{$ i=0 \quad  || \quad e^T\bar{y} - b^T \bar{\lambda} > \epsilon_{opt}$} 
{\label{alg:opt-loop} 
	$j \gets 0$\;
	\While{$j=0 \quad || \quad \lVert vec(dX^*) \rVert_\infty > \epsilon_{apx}$}
	{ \label{alg:apx-loop}
		$\bar{y} \gets \argmin_{}\left\{||y-\bar{y}||: y \in S^P(\bar{X})\right\}$ \; \label{alg:update-y}
		$\bar{\lambda} \gets \argmin\left\{||\lambda -\bar{\lambda}||:\lambda \in S^D(\bar{X})\right\}$ \;\label{alg:update-lambda}
            $dX^*,y^*, \lambda^* \gets \argmin\left\{||vec(dX)||:(dX,y,\lambda) \in S^L(\bar{X},\bar{y},\bar{\lambda})\right\}$	\;
		$\bar{X} \gets \bar{X} + dX^* $\;
		$j\gets j+1$\;
	}
	$\mu \gets 2 \mu $\;
	$i \gets i+1$\;	
}
\end{algorithm}

Since the bilinear terms are omitted in the linearization of both the primal and dual constraints, the solution to the linear problem (\ref{U-model-fixed-penalty}) may not be feasible for the original problem (\ref{U-model-penalty}).
Nonetheless, the step size $dX^*$ tends to increase as the penalty parameter $\mu$ is raised, allowing indirect control over $dX^*$ through adjustments to $\mu$.
If the solution to the linear problem is not feasible for the original problem, then slowing the rate of increase of the penalty parameter  $\mu$ can help reduce the step size $dX^*$ and improve the likelihood of achieving a feasible solution.

\section{A Minimal Example}\label{section:example}
In this section a minimal numerical example is presented to show how the algorithm works.
The notation used here aligns with that in Section~\ref{section:formulation} except the upper-level variables which is denoted by $x$ instead of a matrix $X$.
The upper-level problem is as follows: 
\begin{subequations}\label{example-upper}
\begin{align}
  \min_{x,y_1,y_2}\quad  & |x|  \\
s.t. \quad & y_2 \leq 1.5\\
& y_1, y_2 \in S(x)
\end{align}
\end{subequations}

The absolute value function can easily be replaced with an auxiliary variable and two constraints.
The lower-level problem is as follows:
\begin{subequations}\label{example-lower}
\begin{align}
  \min_{y_1,y_2}\quad  & y_1 + y_2  \\
s.t. \quad & 0.5y_1 + y_2 + x y_1 \geq 3\\
& y_1 + 0.5 y_2 - x y_1 \geq 3 \\
& y_1, y_2 \geq 0
\end{align}
\end{subequations}

The problem is solved to the global optimum using a non-convex QCQP solver.
The optimal solution is $y^*=(2.5, 1.5)$ and $x^*=0.1$.

\begin{figure}[h]
  \centering
  \includegraphics[width=\textwidth]{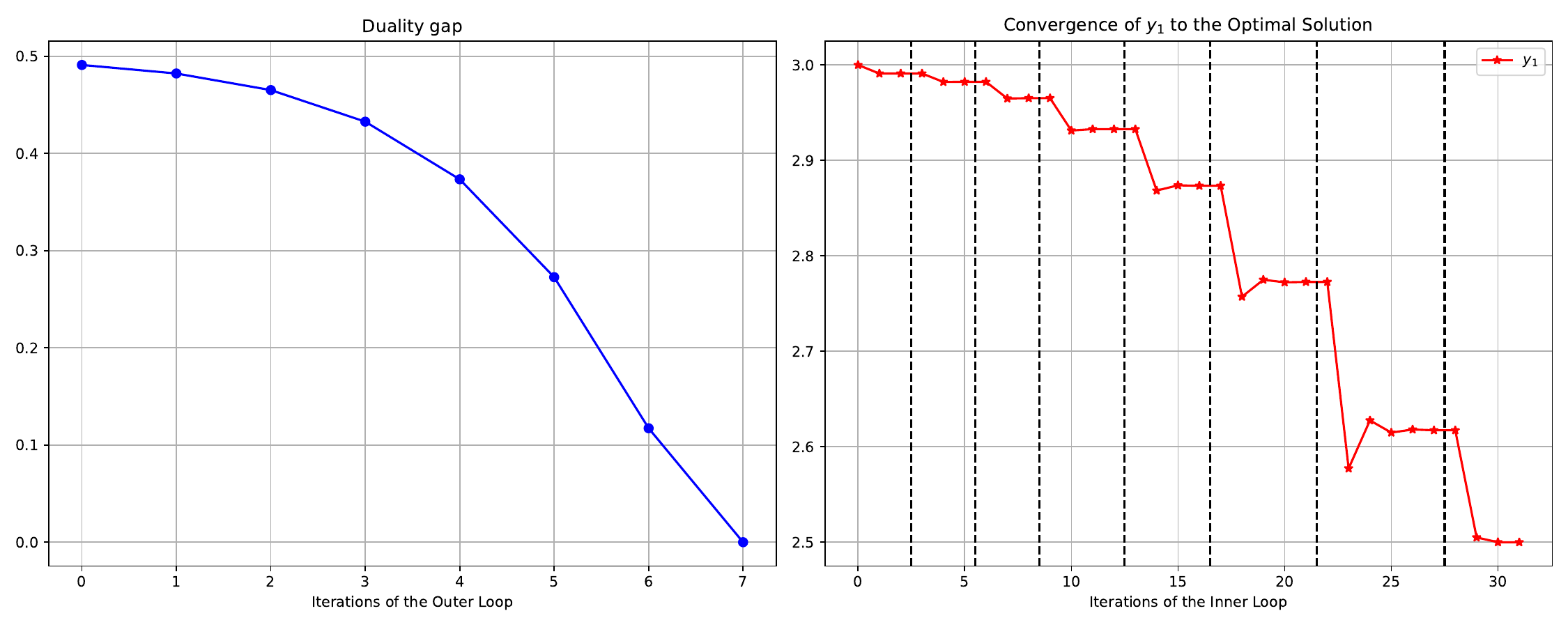}
  \caption{
  (Left)The duality gap of the lower-level problem converges to zero.
  (Right) $\bar{y_1}$ converges to the optimal solution $2.5$. The values are recorded at each iteration of the inner loop and the vertical lines seperates the iterations of the outer loop.
  }
  \label{fig:example}
\end{figure}

Then the Penalty Adaptive Linearization Method is applied to the problem.
The initial value of the upper-level variable is set to zero and the optimal solution of the lower-level problem is $y^*=(2.0, 2.0)$ for this inital value.
The algorithm successfuly converges to the optimal solution of the problem, which is $y^*=(2.5, 1.5)$ and $x^*=0.1$.
Figure \ref{fig:example} shows the convergence of the duality gap of the lower-level problem to zero which means that the algorithm converges to a feasible solution of the main problem.
As expected, the optimal value $\bar{y_2}$ remain constant at $1.5$ throughout the iterations.

\section{Conlcusion}\label{section:conclusion}
In this paper we designed a new algorithm for solving bilevel programming problems with product of upper and lower-level variables in the lower-level problem.
We showed that the algorithm was able to iteratively approximate the bilinear terms in the sample problem and solve it to the global optimum.

The most crucial aspect missing in this algorithm is to prove that the inner loop is guaranteed to converge.
It’s possible that convergence only occurs under specific conditions.
In such cases, identifying non-convergent scenarios can help clarify the algorithm’s limitations, highlight areas for improvement, and specify the conditions under which convergence to a partial minimum is guaranteed.

From a practical perspective, exploring the algorithm’s performance in solving a wide range of real-world problems is of great interest.
However, unlike LP-LP bilevel programming problems, using KKT conditions or strong duality conditions along with the big-M method does not linearize the problem, making it challenging to determine the global optimal solution.

\section*{Acknowledgements}
This research was funded by the German Federal Ministry of Education and Research Infrastructure in project RODES (grant number 05M22RDA). 

\printbibliography

\end{document}